\magnification\magstep1
\baselineskip14pt

\def\Ha {1} 
\def\Hb {2} 
\def\KN {3} 

\def\bmit{\fam=9\relax} 
\font\bmitten=cmmib10
\font\bmitseven=cmmib10 at 7pt
\textfont9=\bmitten \scriptfont9=\bmitseven

\def\bib{\par\noindent\hangindent 25pt}
\def\ldt{\mathrel{.\,.}}
\def\pfbox
  {\hbox{\hskip 3pt\lower2pt\vbox{\hrule
  \hbox to 5pt{\vrule height 7pt\hfill\vrule}
  \hrule}}\hskip3pt}
\def\lf{\lfloor}
\def\rf{\rfloor}
\def\proof{\noindent {\bf Proof.}\enspace}

\centerline{\bf Polynomials involving the floor function}
\centerline{Inger Johanne H{\aa}land\footnote{$^1$}{Agder
College of Engineering, N--4890 Grimstad, Norway} and Donald E. 
Knuth\footnote{$^2$}{Computer Science Department, 
Stanford University, Stanford CA
94305 USA}}
\bigskip

{\narrower\smallskip\noindent
{\bf Abstract.}\enspace
Some identities are presented that generalize the formula
$$x^3=3x\bigl\lf x\lf x\rf\bigr\rf -3\lf x\rf\,\bigl\lf x\lf x\rf\bigr\rf
+\lf x\rf^3+3\{x\}\,\{x\lf x\rf\}+\{x\}^3$$
to a representation of the product $x_0x_1\ldots x_{n-1}$.
\smallskip}

\bigskip\noindent
{\bf 1. Introduction.}\enspace
Let $\lf x\rf$ be the greatest integer less than or equal to~$x$, and let
$\{x\}=x-\lf x\rf$ be the fractional part 
of~$x$. The purpose of this note is to
show how the formulas
$$xy=\lf x\rf y + x\lf y\rf - \lf x\rf\,\lf y\rf+\{x\}\,\{y\}\eqno(1.1)$$
and
$$\eqalignno{xyz&=x\bigl\lf y\lf z\rf\bigr\rf+y\bigl\lf z\lf x\rf\bigr\rf
+z\bigl\lf x\lf y\rf\bigr\rf\cr
&\quad\null -\lf x\rf\,\bigl\lf y\lf z\rf\bigr\rf-\lf y\rf\,\bigl\lf 
z\lf x\rf\bigr\rf-z\bigl\lf x\lf y\rf\bigr\rf\cr
&\quad\null +\lf x\rf\,\lf y\rf\,\lf z\rf\cr
&\quad\null +\{x\}\,\{y\lf z\rf\}+\{y\}\,\{z\lf x\rf\}+\{z\}\,\{x\lf y\rf\}\cr
&\quad\null +\{x\}\,\{y\}\,\{z\}&(1.2)\cr}$$
can be extended to higher-order products $x_0x_1\ldots x_{n-1}$.

These identities make it possible to answer questions about the distribution
mod~1 of sequences having the form
$$\alpha_1n\bigl\lfloor  \alpha_2n \ldots\,\bigl\lfloor\alpha_{k-1}n\,\lfloor
\alpha_kn\rfloor\,\bigr\rfloor\,\ldots\,\bigr\rfloor\,,
\qquad n=1,2,\ldots\,.\eqno(1.3)$$ 
Such sequences are known to be uniformly distributed mod~1 if the real
numbers $1,\alpha_1,\break
\ldots,\alpha_k$ are rationally independent [\Ha]; we will
prove that (1.3) is uniformly distributed in the special case
$\alpha_1=\alpha_2=\cdots =\alpha_k=\alpha$ if and only if $\alpha^k$ is
irrational, when $k$ is prime.
 (It is interesting to compare this result to analogous properties
of the sequence
$$\alpha_0\lfloor \alpha_1n\rfloor\,\lfloor \alpha_2n\rfloor\,\ldots\,
\lfloor\alpha_k n\rfloor\,,\qquad n=1,2\ldots,\eqno(1.4)$$
where $\alpha_0,\alpha_1,\ldots,\alpha_k$ are positive real numbers. If $k\geq
3$, such sequences are uniformly distributed mod~1 if and only if
$\alpha_0$ is irrational~[\Hb].)

\medskip\noindent
{\bf 2. Formulas for the product ${\bmit x_{\bf0} x_{\bf1}\ldots
x_{n-\bf1}}$.}\enspace
The general expression we will derive for $x_0 x_1\ldots x_{n-1}$ contains
$2^{n+1}-n-2$ terms. Given a sequence
$X=(x_0,x_1,\ldots,x_{n-1})$ we regard $x_{n+j}$ as equivalent to~$x_j$, and
for integers $a\leq b$  we define
$$X^{a:b}=\cases{1\,,&if $a=b$;\cr
\noalign{\smallskip}
x_a\lf X^{(a+1):b}\rf\,,&otherwise.\cr}\eqno(2.1)$$
Thus $X^{1:4}=x_1\bigl\lf x_2\lf x_3\rf\bigr\rf$ and
$X^{4:(n+1)}=x_4\bigl\lf x_5\bigl\lf \,\ldots\,\bigl\lf x_{n-1}\lf
x_0\rf\,\bigr\rf \,\ldots\,\bigr\rf\,\bigr\rf$. Using this notation, we obtain
an expression for $x_0 x_1\ldots x_{n-1}$ by taking the sum of
$$\{X^{s_1:s_2}\}\,\{X^{s_2:s_3}\}\,\ldots\,\{X^{s_k:(s_1+n)}\}
-(-1)^k\lf X^{s_1:s_2}\rf\,\lf X^{s_2:s_3}\rf\,\ldots\,\lf X^{s_k:(s_1+n)}\rf
\eqno(2.2)$$
over all nonempty subsets $S=\{s_1,\ldots,s_k\}$ of $\{0,1,\ldots,n-1\}$, where
$s_1<\cdots <s_k$.
This rule defines $2^{n+1}-2$ terms, but in the special case $k=1$ the two
terms of~(2.2) reduce~to
$$\{X^{s_1:(s_1+n)}\}+\lf X^{s_1:(s_1+n)}\rf =X^{s_1:(s_1+n)}\eqno(2.3)$$
so we can combine them and make the overall formula $n$ terms shorter. The
right-hand side of~(1.2) illustrates this construction when $n=3$.

To prove that the sum of all terms (2.2) equals $x_0x_1\ldots x_{n-1}$, we
replace $\{X^{a:b}\}$ by $X^{a:b}-\lf X^{a:b}\rf$ and expand all products. One
of the terms in this expansion is $x_0x_1\ldots x_{n-1}$; it~arises only from
the set $S=\{0,1,\ldots,n-1\}$. The other terms all contain at least one
occurrence of the floor operator, and they can be written
$$x_{u_1}\ldots x_{v_1-1}\lf X^{v_1:u_2}\rf\,x_{u_2}\ldots x_{v_2-1}\,\lf
X^{v_2:u_3}\rf\, x_{u_3}\ldots x_{v_3-1}\;\cdots\,\lf X^{v_k:(u_1+n)}\rf
\eqno(2.4)$$
where $u_1\le v_1<u_2\le v_2<u_3\le \cdots \le v_k<n$. We want to show that all
such terms cancel out. For example, some of the terms in the expansion when
$n=9$ have the form
$$x_1\lf X^{2:4}\rf\,x_4x_5\,\lf X^{6:7}\rf\,\lf X^{7:10}\rf
=x_1\bigl\lf x_2\lf x_3\rf\,\bigr\rf\,x_4x_5\,\lf x_6\rf\,
\bigl\lf x_7\bigl\lf x_8\lf x_0\rf\,\bigr\rf\,\bigr\rf\,,$$
which is (2.4) with $u_1=1$, $v_1=2$, $u_2=4$, $v_2=6$, $u_3=v_3=7$. It is easy
to see that this term arises from the expansion of~(2.2) 
only when $S$ is one of the sets
$\{1,2,4,5,6,7\}$, $\{1,4,5,6,7\}$, $\{1,2,4,5,7\}$, $\{1,4,5,7\}$; in those
cases it occurs with the respective signs $-$, $+$, $+$, $-$, so it does indeed
cancel out.

In general, the only sets $S$ leading to the term (2.4) 
have $S=\{\,s\mid u_j\le
s<v_j\,\}\cup \{\,v_j\mid u_j=v_j\}\cup T$, where $T$ is a subset of
$U=\{\,v_j\mid u_j\neq v_j\,\}$. If $U$ is empty, all parts of the term (2.4)
appear inside floor brackets and this term is cancelled by the second term
of~(2.2). If $U$ contains $m>0$ elements, the $2^m$~choices for~$S$ produce
$2^{m-1}$ terms with a coefficient of $+1$ and $2^{m-1}$ with a coefficient
of~$-1$.
This completes the proof.

Notice that we used no special properties of the floor function in this
argument. The same identity holds when $\lf x\rf$ is an arbitrary function, if
we define $\{x\}=x-\lf x\rf$. 

The formulas become simpler, of course, when all $x_j$ are equal. Let
$$x^{:k}=\cases{1\,,&if $k=0$;\cr
\noalign{\smallskip}
x\lf x^{:(k-1)}\rf\,,&if $k>0$;\cr}\eqno(2.5)$$
and let
$$a_k=\{x^{:k}\}\,,\qquad b_k=\lf x^{:k}\rf\,.\eqno(2.6)$$
Then an identity for $x^n$ can be read off from the coefficients of~$z^n$ in
the formula
$${xz\over 1-xz}\;=\;{a_1z+2a_2z^2+3a_3z^3+\cdots\,\over
1-a_1z-a_2z^2-a_3z^3-\cdots}+{b_1z+2b_2z^2+3b_3z^3+\cdots\,\over
1+b_1z+b_2z^2+b_3z^3+\cdots}\;,\eqno(2.7)$$
which can be derived from (2.2) or proved independently as shown below.
For example,
$$\eqalign{x^2&=a_1^2+2a_2-b_1^2+2b_2\,;\cr
\noalign{\smallskip}
x^3&=a_1^3+3a_1a_2+3a_3+b_1^3-3b_1b_2+3b_3\,;\cr
\noalign{\smallskip}
x^4&=a_1^4+4a_1^2a_2+4a_1a_3+2a_2^2+4a_4\cr
\noalign{\smallskip}
&\qquad\null -b_1^4+4b_1^2b_2-4b_1b_3-2b_2^2+4b_4\,.\cr}$$
In general we have
$$x^n=p_n(a_1,a_2,\ldots,a_n)-p_n(-b_1,-b_2,\ldots,-b_n)\,,\eqno(2.8)$$
where the polynomial
$$p_n(a_1,a_2,\ldots,a_n)=\sum_{k_1+2k_2+\cdots +nk_n=n}\;
{(k_1+k_2+\cdots +k_n-1)!\,n\over k_1!\,k_2!\,\ldots\,k_n!}\;
a_1^{k_1}a_2^{k_2}\,\ldots\,a_n^{k_n}\eqno(2.9)$$
contains one term for each partition of $n$.

It is interesting to note that (2.7) can be written
$${zd\over dz}\,\ln\,{1\over 1-xz}={zd\over dz}\,\ln\,{1\over
1-a_1z-a_2z^2-\cdots \,}-{zd\over dz}\,\ln\,{1\over
1+b_1z+b_2z^2+\cdots\,}\;,$$ 
hence we obtain the equivalent identity
$${1\over 1-xz}={1+b_1z+b_2z^2+b_3z^3+\cdots\,\over
1-a_1z-a_2z^2-a_3z^3-\cdots\,}\;.\eqno(2.10)$$
This identity is easily proved directly, because it says that
$a_k+b_k=xb_{k-1}$ for $k\ge 1$. Therefore it provides an alternative proof of
(2.7). It also 
yields formulas for $x^n$ with mixed $a$'s and $b$'s, and with no negative
coefficients. For example,
$$\eqalign{x^2&=a_1^2+a_2+a_1b_1+b_2\,;\cr
\noalign{\smallskip}
x^3&=a_1^3+2a_1a_2+a_3+(a_1^2+a_2)b_1+a_1b_2+b_3\,;\cr
\noalign{\smallskip}
x^4&=a_1^4+3a_1^2a_2+2a_1a_3+a_2^2+a_4+(a_1^3+2a_1a_2+a_3)b_1\cr
\noalign{\smallskip}
&\qquad\null +(a_1^2+a_2)b_2+a_1b_3+b_4\,.\cr}$$

\medskip\noindent
{\bf 3. Application to uniform distribution.} \enspace
We can now apply the identities to a problem in number theory, as stated in the
introduction. Let $[0\ldt 1)=\{\,x\,\vert\,0\leq x<1\,\}$.

\proclaim
Lemma 1. For all positive integers $k$ and~$l$, there is a function
$f_{k,l}(y_1,y_2,\ldots,y_{k-1})$ from $[0\ldt 1)^{k-1}$ to $[0\ldt1)$ such that
$${x^{:k}\over l}\;\equiv\;
{x^k\over kl}-f_{k,l}\,\left(\left\{{x\over k!\,l}\right\}\,,\,
\left\{{x^2\over k!\,l}\right\}\,,\,\ldots\,,\,\left\{{x^{k-1}\over k!\,l}
\right\}\right)\;\pmod 1\,.\eqno(3.1)$$

\proof
Let
$$\hat{p}_n(a_1,a_2,\ldots,a_{n-1})=p_n(a_1,a_2,\ldots,a_n)-n\,a_n\eqno(3.2)$$
be the polynomial of (2.9) without its (unique) linear term. Then
$${x^{:k}\over l}\,=\,{x^k\over kl}-{1\over kl}\;\hat{p}_k(a_1,\ldots,a_{k-1}) 
+{1\over kl}\;\hat{p}_k(-b_1,\ldots,-b_{k-1})\,.\eqno(3.3)$$
We proceed by induction on $k$, defining the constant $f_{1,l}=0$ for all~$l$.
Then if $y_j=\{x^j\!/k!\,l\}$ and $l_j=k!\,l/j!$ we have
$$a_j=\left\{l_j\;{x^{:j}\over l_j}\right\}\;=\;
\left\{l_j\bigl((j-1)!\,y_j-f_{j,l_j}(y_1,\ldots,y_{j-1})\bigr)\right\}$$
and
$$\eqalign{b_j=\left\lfloor l_j\;{x^{:j}\over l_j}\right\rfloor
&=l_j\left\lfloor{x^{:j}\over l_j}\right\rfloor +\sum_{i=1}^{l_j-1}\,
\left\lfloor\left\{{x^{:j}\over l_j}\right\}+{i\over l_j}\right\rfloor\cr
\noalign{\smallskip}
&\equiv\sum_{i=1}^{l_j-1}\,\left\lfloor\left\{(j-1)!\,y_j-f_{j,l_j}
(y_1,\ldots,y_{j-1})\right\}+{i\over l_j}\right\rfloor\,\pmod{kl}\,,\cr}$$
because of the well-known identities
$$\{lx\}=\bigl\{l\{x\}\bigr\}\,,\qquad \lfloor lx\rfloor=\sum_{i=0}^{l-1}\,
\lfloor x+i/l\rfloor\,,\eqno(3.4)$$
when $l$ is a positive  integer. 
Therefore (3.1) holds with
$$f_{k,l}(y_1,\ldots,y_{k-1})=\left\{{1\over kl}\;\hat{p}_k
(\bar{a}_{1,k,l},\ldots,\bar{a}_{k-1,k,l})-{1\over kl}\;\hat{p}_k
(-\bar{b}_{1,k,l},\ldots,-\bar{b}_{k-1,k,l})\right\}\,,\eqno(3.5)$$
where
$$\eqalignno{\bar{a}_{j,k,l}&=\bigl\{\bigl((j-1)!\,y_j-f_{j,k!l/j!}(y_1,
\ldots,y_{j-1})\bigr)k!\,l/j!\bigr\}\,,&(3.6)\cr
\noalign{\smallskip}
\bar{b}_{j,k,l}&=\sum_{i=1}^{k!l/j!-1}\left\lfloor\{(j-1)!\,y_j-
f_{j,k!l/j!}(y_1,\ldots,y_{j-1})\}+{j!\,i\over k!\,l}\right\rfloor\,.&(3.7)
\cr}$$
For example,
$$\eqalign{f_{2,3}(y)&=\{(\alpha_1^2-\beta_1^2)/6\}\,,\cr
\noalign{\smallskip}
f_{3,1}(y,z)&=\{(3\alpha_1\alpha_2+\alpha_1^3-3\beta_1\beta_2+\beta_1^3)/3\}\,,
\cr}$$
where $\alpha_1=\{6y\}$, $\alpha_2=\{3z-3f_{2,3}(y)\}$, $\beta_1=
\lfloor y+{1\over 6}\rfloor +\lfloor y+{2\over 6}\rfloor +
\cdots +\lfloor y+{5\over 6}\rfloor$, and
$\beta_2=\lfloor\{z-f_{2,3}(y)\}+{1\over 3}\rfloor +\lfloor\{z-f_{2,3}(y)\}
+{2\over 3}\rfloor$. \ \pfbox

\proclaim
Lemma 2. The function $f_{k,l}$ of Lemma 1 does not preserve Lebesgue measure,
and neither does $\{klmf_{k,l}\}$ for any positive integer~$m$.

\proof
It suffices to prove the second statement, for if $f_{k,l}$ were
measure-preserving the functions $\{mf_{k,l}\}$ would preserve Lebesgue measure
for all positive integers~$m$. Notice that
$\{klmf_{k,l}\}=\{m\,\hat{p}_k(\bar{a}_{1,k,l},\ldots,\bar{a}_{k-1,k,l})\}$,
because $\hat{p}_k(-\bar{b}_{1,k,l},\ldots,-\bar{b}_{k-1,k,l})$ is an integer.
The triangular construction of (3.6) makes it clear that
$\bar{a}_{1,k,l},\ldots,\bar{a}_{k-1,k,l}$ are independent random variables
defined on the probability space $[0\ldt 1)^{k-1}$, each uniformly distributed
in $[0\ldt 1)$. Therefore it suffices to prove that
$\{m\,\hat{p}_k(a_1,\ldots,a_{k-1})\}$ is not uniformly distributed when
$a_1,\ldots,a_{k-1}$ are independent uniform deviates.

We can express $\hat{p}_k(a_1,\ldots,a_{k-1})$ in the form
$$\textstyle{k\,a_1a_{k-1}+a_1q_1(a_1,\ldots,a_{k-2})+
k\,a_2a_{k-2}+a_2q_2(a_2,\ldots,
a_{k-3})+\cdots +{1\over 2}k\,a^2_{k/2}\,,}$$
 for some polynomials $q_1,\ldots,
q_{\lfloor (k-1)/2\rfloor}$, where the final term ${1\over 2}k\,a^2_{k/2}$ is
absent when $k$ is odd. Then we can let $y_j=a_j$ for $j\leq {1\over 2}k$ and
$y_j=a_j-q_{k-j}(a_{k-j},\ldots,a_{j-1})/k$ for $j>{1\over 2}k$, obtaining
independent uniform deviates $y_1,\ldots,y_{k-2}$ for which
$m\,\hat{p}_k(a_1,\ldots, a_{k-1})$ equals
$$g_k(y_1,\ldots,y_{k-1})=mk\,y_1y_{k-1}+mk\,y_2y_{k-2}+\cdots 
+({\textstyle{m\over 2}}k\,y^2_{k/2}[k\hbox{ even}])\,.\eqno(3.8)$$
For example, $g_4(y_1,y_2,y_3)=4y_1y_3+2y_2^2$ and
$g_5(y_1,y_2,y_3,y_4)=5y_1y_4+5y_2y_3$ when $m=1$.

The individual terms of (3.8) are independent, and they have monotone
decreasing density functions mod~1.  $\bigl($The density function for the
probability that $\{kxy\}\in[t\ldt t+dt]$ is $\sum_{j=0}^{k-1}\,{1\over
k}\,\ln {k\over j+t}\;dt$.$\bigr)$
 Therefore they cannot
possibly yield a uniform distribution. For if $f(x)$ is the density function
for a random variable on $[0\ldt 1)$, we have $E(e^{2\pi iX})=\int_0^1e^{2\pi
ix}f(x)\,dx\neq 0$ when $f(x)$ is monotone; for example, if $f(x)$ is
decreasing, the imaginary part is $\int_0^{1/2}\sin(2\pi
x)\bigl(f(x)-f(1-x)\bigr)\,dx>0$. If $Y$ is an independent random variable with
monotone density, we have $E(e^{2\pi i\{X+Y\}})=
E(e^{2\pi i(X+Y)})=
E(e^{2\pi iX})E(e^{2\pi
iY})\neq 0$. But $E(e^{2\pi iU})=0$ when $U$ is a uniform deviate. Therefore
(3.8) cannot be uniform mod~1. \ \pfbox

\medskip
Now we can deduce properties of sequences like
$$(\alpha n)^{:k}=\alpha n\bigl\lfloor\alpha n\bigl\lfloor \,\ldots\, \lfloor
\alpha n\rfloor\,\ldots\,\bigr\rfloor\,\bigr\rfloor$$
as $n$ runs through integer values.

\proclaim
Theorem. If the powers $\alpha^2,\ldots,\alpha^{k-1}$ are irrational, the
 sequence $\{m(\alpha n)^k-km(\alpha n)^{:k}\}$, for
$n=1,2,\ldots\,$, is not uniformly distributed in $[0\ldt 1)$ for any
 integer~$m$.

\proof
This result is trivial when $k=1$ and obvious when $k=2$, since $\{(\alpha
n)^2-2(\alpha n)^{:2}\}=\{\alpha n\}^2$. But for large values of~$k$ it seems
to require a careful analysis. By Lemma~1 we have
$$\{m(\alpha n)^k-km(\alpha n)^{:k}\}=\left\{kmf_{k,1}\left(\left\{{\alpha
n\over k!}\right\}\,,\ldots,\,\left\{{\alpha^{k-1}n^{k-1}\over
k!}\right\}\right)\right\}\,,\eqno(3.9)$$ 
and Lemma 2 tells that $\{kmf_{k,1}\}$ is not measure preserving.

Let $S$ be an interval of $[0\ldt 1)$, and $T$ its inverse image in $[0\ldt
1)^{k-1}$ under $\{kf_{k,1}\}$, where $\mu(T)\neq\mu(S)$. 
It is easy to see that if $(y_1,\ldots,y_{k-1})\in T$ and $y_1,\ldots,y_{k-1}$
are irrational, there are values $\epsilon_1,\ldots,\epsilon_{k-1}$ such that
$[y_1\ldt y_1+\epsilon_1)\times\cdots\times[y_{k-1}\ldt
y_{k-1}+\epsilon_{k-1})\subseteq T$. 
Therefore the irrational points of~$T$ can be
covered by disjoint half-open hyperrectangles.
We will show that (3.9) is not uniform by using Theorem 6.4 of [\KN], which
implies that the sequence $(\{\alpha_1n^{e_1}\},\ldots,\{\alpha_sn^{e_s}\})$ is
uniformly distributed in $[0\ldt 1)^s$ whenever $\alpha_1,\ldots,\alpha_s$ are
irrational numbers and the integer exponents $e_1,\ldots,e_s$ are distinct.
Thus
the probability that $\{(\alpha n)^k-k(\alpha n)^{:k}\}\in S$ approaches
$\mu(T)$ as $n\rightarrow\infty$; the distribution is nonuniform. \ \pfbox

\proclaim
Corollary. 
If the powers $\alpha^2,\ldots,\alpha^{k-1}$ are irrational, the
sequence $\{(\alpha n)^{:k}\}$, for $n=1,2,\ldots\,$, is
uniformly distributed in $[0\ldt 1)$ if and only if $\alpha^k$ is irrational.

\proof
If $\alpha^k$ is irrational, $\{\alpha^kn^k\!/k\}$ is uniformly distributed in
$[0\ldt 1)$ and independent of $(\{\alpha
n/k!\},\ldots,\{\alpha^{k-1}n^{k-1}\!/k!\})$, by the theorem quoted above from
[\KN]. Therefore the right-hand side of (3.1) is uniform.

If $\alpha^k$ is rational, say $\alpha^k=p/q$, assume that $\{(\alpha
n)^{:k}\}$ is uniform. Then $\{q(\alpha^kn^k-k(\alpha n)^{:k})\}=\{-qk(\alpha
n)^{:k}\}$ is also uniform, contradicting what we proved. \ \pfbox

\medskip
We conjecture that the theorem and its corollary remain true for all
real~$\alpha$, without the hypothesis that $\alpha^2,\ldots,\alpha^{k-1}$ are
irrational.

\bigskip
\centerline{\bf References}

\medskip
\bib
\Ha
\quad
I. J. H{\aa}land, ``Uniform distribution of generalized polynomials,''
 {\sl Journal of Number Theory\/ \bf 45} (1993), 327--366.

\bib
\Hb
\quad
I. J. H{\aa}land, ``Uniform distribution of generalized polynomials of the
product type,'' {\sl Acta Arithmetica\/ \bf67} (1994), 13--27.

\bib
\KN
\quad
L. Kuipers and H. Niederreiter, {\sl Uniform Distribution of Sequences\/}
(New York: Wiley, 1974).

\bye